# Stability analysis of an SIR model with saturated infection rate and saturated treatment rate

YICHUAN ZHANG [1], CHANGRONG ZHU*
[1]Mathematics and Statistics School, Chongqing University, Chongqing, 401331, P.R. China
*Mathematics and Statistics School, Chongqing University, Chongqing, 401331, P.R. China

**ABSTRACT**: Treatment rate has always been one of the most important factors affecting the spread of infectious diseases. In this paper, we study a treatment function SIR model with treatment rate related to the maximum treatment capacity, whose infection rate is also saturated nonlinear. By calculating and analyzing based on planar dynamic system, the existence and stability of disease-free equilibrium point of the model are found. And it is also found that the situation is more complex for the endemic equilibrium point of the model. Under the conditions of meeting different parameters, the number of endemic equilibrium points may be zero, one or two, and the behavior of the endemic equilibrium point may be a saddle point, a stable or unstable anti-saddle point or center, or even a degenerate case.

## 0 Introduction

Infection rate and treatment rate have always been concerned parameters in infectious disease models. In fact, many scholars are studying the models of contact infection rate and treatment rate. In 1978 , Capasso and Serio [1] introduced the concept of saturation exposure rate $g\left(I\right)S$ [2] in the infectious disease model, where $g\left(I\right)$ is a function of the infection rate associated with $I$ . Liu [3][4] et al. gave the specific form $g\left(I\right)=\dfrac{kSI^{p}}{1+\alpha I^{q}}$ of $g\left(I\right)$ ,where the parameters satisfy $p>0,q>0,\alpha\geq 0$ . E. Amirtharaj [5] studied a class of SIR models considering treatment rates with a transmission rate of $\dfrac{\lambda_1 SI}{1+\alpha_1 I}$ . Since then numerous scholars have investigated the infection rate proposed by Liu with different $p,q$ to analyze the stability of the equilibrium point corresponding to the model. In this process, some scholars have added the treatment rate and made corrections. The most basic treatment function is that the treatment rate is proportional to the number of infected people, which is reasonable when the medical resources are much greater than the needs of the infected people. However, for some explosive infectious diseases, this is not in line with the actual situation. For example, SARS in 2003 instantly infected tens of millions of people. The domestic medical resources at that time were far less than the demand of infected people. Later, W. Wang [6]proposed a segmentation function that gives a maximum treatment rate, and when the treatment volume does not reach the maximum treatment rate, a linear treatment rate is used, and when saturation is reached, the treatment rate no longer increases, which is more reasonable than the linear treatment function. Jing' an Cui[7] et al. introduced Verhulst-type function as the treatment function in the form of $h(I)=\frac{cI}{b+I}$, where $b$ is the number of infected people when $h(I)=\frac{c}{2}$.

* Email address: zhuchangrong@126.com

In this paper, we obtain an SIR model with saturated contact infection rate and saturated treatment rate according to the mechanism of a class of infectious diseases and considering the limited medical resources. The model is as following:

$$\begin{cases} \frac{\mathrm{dS}}{\mathrm{dt}} = A - \mathrm{d}S - \frac{\lambda SI}{1+\alpha I} \\ \frac{\mathrm{dI}}{\mathrm{dt}} = \frac{\lambda SI}{1+\alpha I} - (d+m)I - \frac{cI}{b+I} \\ \frac{\mathrm{dR}}{\mathrm{dt}} = \frac{cI}{b+I} - \mathrm{d}R + mI \end{cases} \tag{0.1}$$

where $S$ is susceptible population, $I$ is infected population, $R$ is removed population, $A$ is birth rate, $d$ is death rate, $\lambda$ is proportional constant, $\alpha$ is influencing parameter, $m$ is natural recovery rate, $\frac{cI}{b+I}$ is medical recovery rate of infected people, $c$ is the maximum recovery number per unit time, and $b$ represents the number of infected people when $\frac{cI}{b+I}$ is taken as $\frac{\mathrm{c}}{2}$. Here we assume that all parameters are positive numbers.

To facilitate the study, we first simplify the model. We use $N(t)=I(t)+R(t)+S(t)$ to represent the total population, where $N$ is a function of $t$. Add all the equations of $(0.1)$ and we get $\frac{dN}{dt}=A-dN$. It can be found that for any initial condition, when time $t$ tends to infinity, $N(t)$ will tend to constant $\frac{A}{d}$, so we can discuss the system composed of the first two equations:

$$\begin{cases} \frac{\mathrm{dS}}{\mathrm{dt}} = A - \mathrm{d}S - \frac{\lambda SI}{1+\alpha I} \\ \frac{\mathrm{dI}}{\mathrm{dt}} = \frac{\lambda SI}{1+\alpha I} - (d+m)I - \frac{cI}{b+I} \end{cases} \tag{0.2}$$

Next, we analyze the existence and stability of the equilibrium point of model $(0.2)$ in detail.

# 1 The existence of equilibrium point

To simplify the calculation, we first introduce an equivalent transformation for system

$(0.2)$. Denote $k = d + m$ and let $S = \frac{k}{\lambda} x, I = \frac{k}{\lambda} y, t = \frac{1}{k} \tau$. Then system $(0.2)$ becomes:

$$\begin{cases} \dot{x} = a - ex - \dfrac{xy}{1+gy} \\ \dot{y} = \dfrac{xy}{1+gy} - y - \dfrac{my}{b+ny} \end{cases} \tag{1.1}$$

where $a = \frac{A\lambda}{k^2}, e = \frac{d}{k}, g = \frac{\alpha k}{\lambda}, m = \frac{ck}{\lambda}, n = \frac{k}{\lambda}$ and $0 < e < 1$. The equilibrium points of system $(1.1)$ is the solution of the following equation:

$$\begin{cases} a - ex - \dfrac{xy}{1+gy} = 0 \\ \dfrac{xy}{1+gy} - y - \dfrac{my}{b+ny} = 0 \end{cases} \tag{1.2}$$

Obviously, there is a disease-free equilibrium point $E_0 = \left( \frac{a}{e}, 0 \right)$ for system $(1.1)$ regardless of the value of the parameters. The basic regeneration number of infectious disease model $(1.1)$ calculated by the symbols and methods in the reference [8] was: $R_0 = \frac{a}{e\left(1 + \frac{m}{b}\right)} = \frac{ab}{e(m+b)}$.

Next, consider the case of the endemic equilibrium point. From the first equation of equation $(1.2)$, we get $x = \frac{a(1+gy)}{e + (eg+1)y}$. Then we substitute it into the second equation of equation $(1.2)$ to simplify it and get:

$$a_1 y^2 + a_2 y + a_3 = 0 \tag{1.3}$$

where $a_1 = (eg+1)n > 0, a_2 = (b+m)(eg+1) + n(e-a)$,

$a_3 = e(b+m) - ab = e(b+m)(1 - R_0)$.

Obviously, equation $(1.3)$ is a quadratic equation with respect to $y$. The discriminant is $\Delta = a_2^2 - 4a_1 a_3$.

Next, we discuss the existence of positive solutions of equation $(1.3)$ by situation.

**Lemma 1.1**: For any $g, b, m > 0$ and $0 < e < 1$, then we have

(1) For any $n>0$, when $a>\frac{e(m+b)}{b}$, then the equation $(1.3)$ must have a positive root $y_1=\frac{-a_2+\sqrt{\Delta}}{2a_1}$;

(2) When $n>\frac{(eg+1)b^2+m(eg+1)b}{me}$ and $a=\frac{e(m+b)}{b}$, then the equation $(1.3)$ must have a positive root $y_1=\frac{-a_2+\sqrt{\Delta}}{2a_1}$.

**Proof**: (1) When $a>\frac{e(m+b)}{b}$, there is $R_0>1$, then $a_3<0$, so $\Delta=a_2^2-4a_1a_3>0$. So we know that the equation $(1.3)$ has a positive root, which is called $y_1$. Solve equation $(1.3)$ to get $y_1=\frac{-a_2+\sqrt{\Delta}}{2a_1}$.

(2) When $a=\frac{e(m+b)}{b}$, there is $R_0=1$, then $a_3=0$, so $\Delta=a_2^2>0$. We notice that equation $(1.3)$ has a positive root if and only if $a_2<0$. From $a_2<0$,we get $\frac{(b+m)(eg+1)+ne}{n}<a$, and we have $a=\frac{e(m+b)}{b}$, then the two equations are combined to eliminate $a$ and simplified to obtain:

$$(eg+1)b^2+m(eg+1)b-mne<0 \qquad (1.4)$$

Inequality $(1.4)$ holds if and only if $n>\frac{(eg+1)b^2+m(eg+1)b}{me}$. Therefore, equation $(1.3)$ has positive roots $y_1=\frac{-a_2+\sqrt{\Delta}}{2a_1}$ when $a=\frac{e(m+b)}{b}$ and $n>\frac{(eg+1)b^2+m(eg+1)b}{me}$.

The proof is complete.

Let $a_4=\frac{\left[en-(b-m)(eg+1)\right]+\sqrt{4m(eg+1)\left[ne-b(eg+1)\right]}}{n}$.

**Lemma 1.2**: For any $g,b,m>0$ and $0<e<1$, then we have

(1) When $n>\frac{(eg+1)b^2+m(eg+1)b}{me}$ and $a_4<a<\frac{e(m+b)}{b}$, the equation $(1.3)$ has two distinct positive roots $y_1=\frac{-a_2+\sqrt{\Delta}}{2a_1}$ and $y_2=\frac{-a_2-\sqrt{\Delta}}{2a_1}$.

(2) When $n > \frac{(eg+1)b^2 + m(eg+1)b}{me}$ and $a = a_4$, the equation $(1.3)$ has a positive root $y_3 = -\frac{a_2}{2a_1}$.

(3) When $n > \frac{b(eg+1)}{e}$ and $0 < a < a_4$ or $a_4 < a < \frac{e(m+b)}{b}$ and $a_4 < a < \frac{e(m+b)}{b}$, equation $(1.3)$ has no positive root.

**Proof:** When $a < \frac{e(m+b)}{b}$, there is $R_0 < 1$, then $a_3 > 0$. We notice that equation $(1.3)$ has a positive root if and only if $a_2 < 0$ and $\Delta \geq 0$. It follows from Lemma 1.1 that $a_3 > 0$, $a_2 < 0$ are equivalent to the parameters $a$ and $n$ satisfying $n > \frac{(eg+1)b^2 + m(eg+1)b}{me}$ and

$$\frac{(b+m)(eg+1)+ne}{n} < a < \frac{e(m+b)}{b} \tag{1.5}$$

Since the positive and negative aspects of $\Delta$ are involved, we specify $\Delta$ and obtain the following form by simplification:

$$\Delta = f(a) = n^2a^2 + 2n\left[(b-m)(eg+1) - en\right]a + \left[(b+m)(eg+1) - ne\right]^2,$$

Obviously $f(a)$ is a quadratic function with respect to $a$. The discriminant is

$$\begin{aligned}\Delta_1 &= \left[2n(b-m)(eg+1) - 2en^2\right]^2 - 4n^2\left[(b+m)(eg+1) - ne\right]^2 \\ &= 16\mathrm{n}^2 m(eg+1)\left[ne - b(eg+1)\right]\end{aligned}$$

When $n > \frac{(eg+1)b^2 + m(eg+1)b}{me}$, we have $n > \frac{b(eg+1)}{e}$, then $\Delta_1 > 0$ and $(b-m)(eg+1) - en < 0$, solving for the roots of $f(a)$ is:

$$a_{\pm} = \frac{\left[en - (b-m)(eg+1)\right] \pm \sqrt{4m(eg+1)\left[ne - b(eg+1)\right]}}{n}.$$

Because $a$ has to satisfy the inequalities $(1.5)$ and $n > \frac{(eg+1)b^2 + m(eg+1)b}{me}$, it is known that $a_+$ satisfies the condition and makes $\Delta = f(a) = 0$. Note $a_4 = a_+$. When $a \in (0, a_4)$, there is $\Delta = f(a) < 0$; When $a = a_4$, there is $\Delta = f(a) = 0$; When $a \in \left(a_4, \frac{e(m+b)}{b}\right)$, there is

$\Delta = f(a) > 0$. Next, we illustrate the case of positive roots of the equation $(1.3)$ in three cases:

Case 1: When $n > \frac{(eg+1)b^2 + m(eg+1)b}{me}$ and $a_4 < a < \frac{e(m+b)}{b}$, it is known that $\Delta > 0$, $a_3 > 0$, $a_2 < 0$, which shows that the equation $(1.3)$ has two distinct positive roots $y_1$ and $y_2$, which yields $y_1 = \frac{-a_2 + \sqrt{\Delta}}{2a_1}$ and $y_2 = \frac{-a_2 - \sqrt{\Delta}}{2a_1}$.

Case 2: When $n > \frac{(eg+1)b^2 + m(eg+1)b}{me}$ and $a = a_4$, it is known that $\Delta = 0$, $a_3 > 0$, $a_2 < 0$, which shows that the equation $(1.3)$ has a positive root $y_3 = -\frac{a_2}{2a_1}$.

Case 3: When $n > \frac{b(eg+1)}{e}$ and $0 < a < a_4$ or $n < \frac{b(eg+1)}{e}$ and $a_4 < a < \frac{e(m+b)}{b}$, it is known that $\Delta < 0$, $a_3 > 0$, which shows that the equation $(1.3)$ has no roots, and thus no positive roots. End of proof.

Let $x_1 = \frac{e(1+gy_1)}{e+(eg+1)y_1}$, $x_2 = \frac{e(1+gy_2)}{e+(eg+1)y_2}$, $x_3 = \frac{e(1+gy_3)}{e+(eg+1)y_3}$.

With the above analysis and based on the conclusions of Lemma 1.1 and Lemma 1.2, we arrive at the situation regarding the endemic equilibrium point of system $(1.1)$. The following theorem is concluded:

**Theorem 1.3:** For any $g, b, m > 0$ and $0 < e < 1$, then we have

(1) For any $n > 0$, if $a > \frac{e(m+b)}{b}$, then system $(1.1)$ has an endemic equilibrium point $E_1(x_1, y_1)$;

(2) If $n > \frac{(eg+1)b^2 + m(eg+1)b}{me}$ and $a = \frac{e(m+b)}{b}$, then system $(1.1)$ has an endemic equilibrium point $E_1(x_1, y_1)$;

(3) If $n > \frac{(eg+1)b^2 + m(eg+1)b}{me}$ and $a_4 < a < \frac{e(m+b)}{b}$,then system $(1.1)$ has two endemic equilibrium points $E_1(x_1, y_1)$ and $E_2(x_2, y_2)$;

(4) If $n > \frac{(eg+1)b^2 + m(eg+1)b}{me}$ and $a = a_4$, then system $(1.1)$ has a double endemic equilibrium point $E_3(x_3, y_3)$;

(5) When $n > \frac{b(eg+1)}{e}$ and $0 < a < a_4$ or $a_4 < a < \frac{e(m+b)}{b}$ and $a_4 < a < \frac{e(m+b)}{b}$, there is no endemic equilibrium in system $(1.1)$.

## 2 Stability of equilibrium point

The equilibrium point is of great significance to the study of infectious disease model. Next, we study the stability of disease-free equilibrium point and endemic equilibrium point of the model. First, we study the disease-free equilibrium.

**Theorem 2.1:**(1) If $a < \frac{e(m+b)}{b}$ ,then The disease-free equilibrium point $E_0\left(\frac{a}{e}, 0\right)$ of system $(1.1)$ is locally asymptotically stable and a stable node; If $a > \frac{e(m+b)}{b}$, the disease-free equilibrium point $E_0$ is saddle point;

(2) If $a = \frac{e(m+b)}{b}$ and $n \neq \frac{(eg+1)b^2 + m(eg+1)b}{me}$ ,then disease-free equilibrium point $E_0$ is a saddle-node.

**Proof:** The Jacobi matrix of system $(1.1)$ at $E_0$ is

$$J_{E_0} = \begin{pmatrix} -e & \frac{a}{e} \\ 0 & \frac{a}{e} - 1 - \frac{m}{b} \end{pmatrix}$$

and the corresponding eigenvalue is $\lambda_1 = -e < 0$ and $\lambda_2 = \frac{a}{e} - 1 - \frac{m}{b}$.

(1) According to reference [9], when $a < \frac{e(m+b)}{b}$, so $\lambda_2 = \frac{a}{e} - 1 - \frac{m}{b} < 0$, then $E_0$ is locally asymptotically stable and a stable node; When $a > \frac{e(m+b)}{b}$, so $\lambda_2 > 0$, then $E_0$ is the saddle point.

(2) When $a = \frac{e(m+b)}{b}$,so $\lambda_1 < 0$ and $\lambda_2 = 0$,then the $E_0$ is a Lyapunov singularity. Firstly,

move $E_0$ to the origin and transform $\begin{cases} u = x - \frac{a}{e} \\ v = y \end{cases}$ to get system

$$\begin{cases} \frac{du}{dt} = -eu - \frac{uv}{1+gv} - \frac{\frac{a}{e}v}{1+gv} \\ \frac{dv}{dt} = \frac{uv}{1+gv} + \frac{\frac{a}{e}v}{1+gv} - v - \frac{mv}{b+nv} \end{cases} \qquad (2.1)$$

Taylor expansion of system $(2.1)$ at the origin gives:

$$\begin{cases} \frac{du}{dt} = -eu - \frac{a}{e}v + \varphi(u,v) \\ \frac{dv}{dt} = \phi(u,v) \end{cases} \qquad (2.2)$$

where $\varphi(u,v) = -uv + \frac{ag}{e}v^2 + o\left((u,v)^3\right)$ and $\phi(u,v) = uv + \left(\frac{mn}{b^2} - \frac{ag}{e}\right)v^2 + o\left((u,v)^3\right)$.

From equation $-eu - \frac{a}{e}v + \varphi(u,v) = 0$, using implicit function theorem [10], we get: $u = h(v)$ and $h(0) = 0$ .Let $u = h(v) = a_1 v + a_2 v^2 + o(v^3)$ . Substituting it into the equation and comparing the coefficients, we get:

$$u = h(v) = -\frac{a}{e^2}v + \frac{a(eg+1)}{e^3}v^2 + o(v^3)$$

A new variable $z = u - h(v)$ is introduced, and the formula $(2.2)$ is changed into:

$$\begin{cases} \frac{dz}{dt} = -ez - zv + o\left((z,v)^2\right) \\ \frac{dv}{dt} = z + g_1 v^2 + o\left((z,v)^2\right) \end{cases}$$

Where $g_1 = \frac{mn}{b^2} - \frac{a(eg+1)}{e^2}$. If $a \neq \frac{mne^2}{b^2(eg+1)}$, then $g_1 \neq 0$, so $E_0$ is a saddle node. Since $a = \frac{e(m+b)}{b}$ , $\frac{e(m+b)}{b} \neq \frac{mne^2}{b^2(eg+1)}$ holds if and only if $n \neq \frac{(eg+1)b^2 + m(eg+1)b}{me}$ .Complete the certificate.

Generally speaking, to consider the stability at the equilibrium point, it is necessary to consider the real part of the eigenvalue of Jacobi matrix at the equilibrium point. If $(x, y)$ is the endemic equilibrium point of system $(1.1)$, the Jacobi matrix is

$$J=\begin{pmatrix} -e-\dfrac{y}{1+gy} & -\dfrac{x}{(1+gy)^2} \\ \dfrac{y}{1+gy} & \dfrac{x}{(1+gy)^2}-1-\dfrac{mb}{(b+ny)^2} \end{pmatrix}$$

and its characteristic equation is $\lambda^2-tr(J)\lambda+\det(J)=0$ ,where $tr(J)$ and $\det(J)$ represent the trace and determinant of matrix $J$ respectively, and their expressions are:

$$tr(J)=-e-\frac{y}{1+gy}+\frac{x}{(1+gy)^2}-1-\frac{mb}{(b+ny)^2},$$

$$\det(J)=-\frac{ex}{(1+gy)^2}+e+\frac{mbe}{(b+ny)^2}+\frac{y}{1+gy}+\frac{mby}{(1+gy)(b+ny)^2}.$$

For the endemic equilibrium point $E_i$ , where $i=1,2,3$ , they all satisfy equation $\frac{x_i}{1+gy_i}=1+\frac{m}{b+ny_i}$. Substitute the expression into $tr(J)$ and $\det(J)$ to get:

$$tr(J)=\frac{-H(y_i)}{(1+gy_i)(b+ny_i)^2},\det(J)=\frac{y_iM(y_i)}{(1+gy_i)(b+ny_i)^2},$$

where $H(y_i)=\left[e+(eg+g+1)y_i\right](b+ny_i)^2-m(n-bg)y_i$,

$M(y_i)=(eg+1)(b+ny_i)^2+m\left[b(eg+1)-ne\right]$.

Next, we explain the stability of the equilibrium point of endemic diseases in different situations.

**Theorem 2.2:** If $n>\frac{(eg+1)b^2+m(eg+1)b}{me}$ and $a=a_4$, then the equilibrium point $E_3$ is degenerate.

**Proof:** If $n>\frac{(eg+1)b^2+m(eg+1)b}{me}$ and $a=a_4$, we can calculate $M(y_3)=0$, which means $\det\left(J_{E_3}\right)=0$. So, the equilibrium point $E_3$ is degenerate. The certificate is over.

**Remark:** There are two cases of zero eigenvalue of system $(1.1)$ at equilibrium point $E_3$ .Let

$N=(eg+1)\left[ne-b(eg+1)\right]$.

(1) If $neN+(N-n)\left[\sqrt{mN}-b(eg+1)\right]\neq 0$ , then the system $(1.1)$ has a zero eigenvalue at $E_3$;

(2) If $neN+(N-n)\left[\sqrt{mN}-b(eg+1)\right]=0$ ,then the system $(1.1)$ has double zero eigenvalues at $E_3$.

**Theorem 2.3**: If $n>\dfrac{(eg+1)b^2+m(eg+1)b}{me}$ and $a_4<a<\dfrac{e(m+b)}{b}$ , then the endemic equilibrium point $E_2(x_2,y_2)$ of system $(1.1)$ is saddle point.

**Proof :** If $n>\dfrac{(eg+1)b^2+m(eg+1)b}{me}$ ,since $M'(y)=2(eg+1)n^2y+2bn(eg+1)>0$ , $M(y)$ is monotonically increasing with respect to $y$ .Because of $y_2=\dfrac{-a_2-\sqrt{\Delta}}{2a_1}$ and $y_3=-\dfrac{a_2}{2a_1}$ , we have $y_2<y_3$ . Theorem 2.2 shows that: $M(y_3)=0$ , which implies that $M(y_2)<0$, so $\det\left(J_{E_2}\right)<0$. Therefore, the equilibrium point stability theory of plane system shows that the endemic equilibrium point $E_2(x_2,y_2)$ is a saddle point. The certificate is over.

The stability of equilibrium point $E_1$ is complex, and different situations will appear with different parameters. Before explaining the stability of equilibrium point $E_1$, let $w=s_1\sqrt{\Delta}+s_2$, $l=n(eg+g+1)$, where

$$s_1=\left(a_1a_3-a_2^2\right)nl+a_1a_2\left(en^2+2bl\right)-a_1^2\left(\frac{b^2l}{n}+2ben+bmg-mn\right)$$

$$s_2=\left(a_2^3-3a_1a_2a_3\right)nl+\left(2a_1^2a_3-a_1a_2^2\right)\left(en^2+2bl\right)+a_1^2a_2\left(\frac{b^2l}{n}+2ben+bmg-mn\right)-2a_1^3eb^2$$

**Theorem 2.4**: If $n>\dfrac{(eg+1)b^2+m(eg+1)b}{me}$ and $a_4<a\leq\dfrac{e(m+b)}{b}$ , the stability of

equilibrium point $E_1$ has the following three situations:

(1) If $w<0$, then $E_1$ is a stable node or focus;

(2) If $w=0$, then $E_1$ is the weak center;

(3) If $w>0$, then $E_1$ is an unstable node or focus.

**Proof:** Firstly, we prove that $\det\left(J_{E_1}\right)>0$. According to Theorem 2.3, $M\left(y\right)$ is monotonically increasing with respect to $y$. Because of $y_1=\dfrac{-a_2+\sqrt{\Delta}}{2a_1}$ and $y_3=-\dfrac{a_2}{2a_1}$, we have $y_3<y_1$. Theorem 2.2 shows that: $M\left(y_3\right)=0$, so $M\left(y_1\right)>0$, in other words $\det\left(J_{E_1}\right)>0$.

Next, we study the symbol of $H\left(y_1\right)$. We notice that $y_1$ satisfies the equation $(1.3)$ and $y_1=\dfrac{-a_2+\sqrt{\Delta}}{2a_1}$. Substitute the condition into the expression of $H\left(y_1\right)$, and simplify it to get:

$$H\left(y_1\right)=\left(a_1y_1^2+a_2y_1+a_3\right)\left(b_1y+b_2\right)-\frac{w}{2a_1^3}=-\frac{w}{2a_1^3},$$

where $b_1=\dfrac{nl}{a_1}, b_2=\dfrac{1}{a_1}\left(en^2+2bl-\dfrac{a_2nl}{a_1}\right)$.

Therefore, if $w<0$, we get $H\left(y_1\right)>0$, then $tr\left(J_{E_1}\right)<0$, so $E_1$ is a stable node or focus; if $w=0$, we get $H\left(y_1\right)=0$, then $tr\left(J_{E_1}\right)=0$ ,so $E_1$ is the weak center; if $w>0$, we get $H\left(y_1\right)<0$, then $tr\left(J_{E_1}\right)>0$, so $E_1$ is a unstable node or focus.

**Theorem 2.5:** The equilibrium point $E_1\left(x_1,y_1\right)$ is locally asymptotically stable when the system meets one of the following conditions:

(1) $n\le\dfrac{(eg+1)b^2+m(eg+1)b}{me}$ 且 $a>\dfrac{e(m+b)}{b}$,

(2) $n>\dfrac{(eg+1)b^2+m(eg+1)b}{me}$ 且 $a>\dfrac{mne^2}{b^2(eg+1)}$.

**Proof:** According to the theorem in literature [11], consider the system $\dfrac{dx}{dt}=f\left(x,\phi\right)$ with parameter $\phi$. Where $f$ is at least second order continuously differentiable with respect to $x$ and $\phi$. Now let $\phi = m$ be the bifurcation parameter, such that $R_0 > 1$ when $\phi > \frac{b(a-e)}{e}$, $R_0 < 1$ when $\phi < \frac{b(a-e)}{e}$, and $x_0 = (0,0)$ is the disease-free equilibrium point with respect to any $\phi$. In order to ensure that the disease-free equilibrium point is at $(0,0)$, we transform $x_1 = x - \frac{a}{e}, x_2 = y$, and system $(1.1)$ becomes

$$\begin{cases} \frac{dx_1}{dt} = -ex_1 - \frac{x_1 x_2}{1+gx_2} - \frac{\frac{a}{e}x_2}{1+gx_2} \triangleq f_1 \\ \frac{dx_2}{dt} = \frac{x_1 x_2}{1+gx_2} + \frac{\frac{a}{e}x_2}{1+gx_2} - x_2 - \frac{mx_2}{b+nx_2} \triangleq f_2 \end{cases} \tag{2.3}$$

The Jacobi matrix of the system $(2.3)$ at the parameter $m = \frac{b(a-e)}{e}$ and the disease-free equilibrium $x_0 = (0,0)$ is calculated as follows:

$$J\left(x_0, \frac{b(a-e)}{e}\right) = \begin{pmatrix} -e & -\frac{a}{e} \\ 0 & 0 \end{pmatrix}$$

The eigenvalues of this matrix is $\lambda_1 = -e < 0$ and $\lambda_2 = 0$, then the first condition of the theorem in literature [11] is satisfied.

Now, if $\varpi = (\omega_1, \omega_2)^T$ is defined as the right eigenvector of the eigenvalue $\lambda_2 = 0$, then

$$\begin{pmatrix} -e & -\frac{a}{e} \\ 0 & 0 \end{pmatrix} (\omega_1, \omega_2)^T = 0.$$

Through calculation, it can be obtained that $\varpi = \left(1, -\frac{e^2}{a}\right)^T$.

If $v = (v_1, v_2)$ is defined as the left eigenvector of the eigenvalue $\lambda_2 = 0$, calculate by the same method to get $v = (0,1)$.

Calculate the partial derivative at the disease-free equilibrium point $x_0 = (0,0)$ and get:

$$\frac{\partial^2 f_1}{\partial x_1 \partial x_2} = -1, \frac{\partial^2 f_1}{\partial x_2 \partial x_1} = -1, \frac{\partial^2 f_1}{\partial x_2 \partial x_2} = \frac{2ag}{e},$$

$$\frac{\partial^2 f_2}{\partial x_1 \partial x_2}=1, \frac{\partial^2 f_2}{\partial x_2 \partial x_1}=1, \frac{\partial^2 f_2}{\partial x_2 \partial x_2}=-\frac{2ag}{e}+\frac{2mn}{b^2}, \frac{\partial^2 f_2}{\partial x_2 \partial m}=-\frac{1}{b},$$

The remaining second-order partial derivatives are all zero. According to Theorem 4.1 in reference [11], $\hat{a},\hat{b}$ is defined as:

$$\hat{a}=\sum_{k,i,j=1}^{2} v_k \omega_i \omega_j \frac{\partial^2 f_k}{\partial x_i \partial x_j}(x_0,0), \hat{b}=\sum_{k,i=1}^{2} v_k \omega_i \frac{\partial^2 f_k}{\partial x_i \partial \phi}(x_0,0).$$

Calculated to:

$$\hat{a}=v_2\omega_1\omega_2 \frac{\partial^2 f_2}{\partial x_1 \partial x_2}+v_2\omega_2\omega_1 \frac{\partial^2 f_2}{\partial x_2 \partial x_1}+v_2\omega_2\omega_2 \frac{\partial^2 f_2}{\partial x_2 \partial x_2}=-\frac{2e^2}{a}+\frac{e^4}{a^2}\left(\frac{2mn}{b^2}-\frac{2ag}{e}\right),$$

$$\hat{b}=v_2\omega_2 \frac{\partial^2 f_2}{\partial x_2 \partial m}=\frac{e^2}{ab}>0.$$

When $a>\frac{mne^2}{b^2(eg+1)}$, we have $\hat{a}<0$. It is known from Theorem 4.1 in reference [11] that the only endemic equilibrium point $E_1$ near the disease-free equilibrium $E_0$ is locally asymptotically stable. Then combined with the existence condition of equilibrium point $E_1$: $a>\frac{e(m+b)}{b}$, We get: if $n>\frac{(eg+1)b^2+m(eg+1)b}{me}$, then $a>\frac{mne^2}{b^2(eg+1)}>\frac{e(m+b)}{b}$; if $n\le\frac{(eg+1)b^2+m(eg+1)b}{me}$, then $a>\frac{e(m+b)}{b}\ge\frac{mne^2}{b^2(eg+1)}$. To sum up, when one of the two situations is established, the endemic equilibrium point $E_1$ is locally asymptotically stable. The certificate is over.

## 3 Conclusion

In this paper, a kind of SIR epidemic model with saturated infection rate and saturated treatment rate is studied by using the qualitative theory of plane system in differential dynamic system. Firstly, the parameter conditions of the existence of disease-free equilibrium and endemic equilibrium are obtained through research and analysis. Then the possible stability states of the system corresponding to different equilibrium points with different parameters under different parameter conditions are obtained by combining qualitative theory.

## References

[1] Fred Brauer,Carlos Castillo Chavez. Mathematical Models in Population Biology and

Epidemiology[M] Springer, New York, NY:2011.

[2] Ma Zhien, Zhou Yicang, Wang Wendi, Jin Zhen. Mathematical modeling and research of infectious disease dynamics [M]. Beijing: Science Press, 2004.

[3] Liu W M,Hethcote H W,Levin S A. Dynamical behavior of epidemiological models with nonlinear incidence rates.[J]. Journal of mathematical biology,1987,25(4).

[4] Liu W M,Levin S A,Iwasa Y. Influence of nonlinear incidence rates upon the behavior of SIRS epidemiological models.[J]. Journal of mathematical biology,1986,23(2)

[5] D.Jasmine, E.C.Henry，Amirtharaj, Backward bifurcation of SIR epidemic model with nonmonotonic incidence rate under treatment[J], IOSR-JM. 10 (2014) 22-32.

[6] Shigui Ruan,Wendi Wang. Dynamical behavior of an epidemic model with a nonlinear incidence rate[J]. Journal of Differential Equations,2003,188(1).

[7] Jinǵan Cui,Xiaoxia Mu,Hui Wan. Saturation recovery leads to multiple endemic equilibria and backward bifurcation[J]. Journal of Theoretical Biology,2008,254(2).

[8] P. van den Driessche,James Watmough. Reproduction numbers and sub-threshold endemic equilibria for compartmental models of disease transmission[J]. Mathematical Biosciences,2002,180(1).

[9] John Guckenheimer,Philip Holmes. Nonlinear Oscillations, Dynamical Systems, and Bifurcations of Vector Fields[M] Springer, New York, NY.

[10] Shui Nee Chow,Jack K. Hale. Methods of Bifurcation Theory[M] Springer, New York, NY.

[11] Castillo-Chavez Carlos,Song Baojun. Dynamical models of tuberculosis and their applications .[J]. Mathematical biosciences and engineering : MBE,2004,1(2).